\input amstex
\documentstyle{amsppt}
\TagsOnRight
\pagewidth{5.5 in}
\pageheight{8 in}
\magnification=1200
\NoBlackBoxes
  
%%%%%%%%%%%%%%%%%%%%%%%%%%%%%%%%%%%%%%%%%%%%%  MACROS BEGIN  %%%%%%%%%%%%%%
\def\e{{\roman e}}

\def\conj#1{\overline{#1}}
\def\inner#1#2#3{ \langle #1,\, #2 \rangle_{_{#3}}}
\def\Dinner#1#2{ \langle #1,\, #2 \rangle_{D}}
\def\Dpinner#1#2{ \langle #1,\, #2 \rangle_{D^\perp}}

\def\vartau{\tau\hskip2pt\!\!\!\!\iota\hskip1pt}
%\def\vartau{\boldsymbol\tau\hskip2pt\!\!\!\!\boldsymbol\iota\hskip1pt}

%%%%%%%%%%%%%%%%  References Macros:
\def\FB{1}
\def\EE{2}
\def\EL{3}
\def\MRa{4}
\def\MRb{5}
\def\SWa{6}
\def\SWb{7}
\def\SWc{8}
\def\SWd{9}
%%%%%%%%%%%%%%%%%%%%%%%%%%%%%%%%%%%%%%%%%%%%%  MACROS END  %%%%%%%%%%%%%%

%%%%%%%%%%%%%%%%%%%%%%%%%%%%%%%%%%%%%%%%%%%%%%%%%%%%%%%%%%%%%%%%%%%%%%%%%
%%
%%  This copy submitted to the Mathematische Annalen, Dr. Wolfgang Lueck,
%%  on Wed Dec 6, 2000.
%%
%%%%%%%%%%%%%%%%%%%%%%%%%%%%%%%%%%%%%%%%%%%%%%%%%%%%%%%%%%%%%%%%%%%%%%%%%

\topmatter

\title
On Fourier orthogonal projections in the rotation algebra 
\endtitle

%\rightheadtext{Fourier orthogonal projections}

\author
S.~Walters
\endauthor

\affil
\sevenrm The University of Northern British Columbia
\endaffil

\date
\sevenrm November 30, 2000
\enddate

\address
Department of Mathematics and Computer Science, The University
of Northern British Columbia, Prince George, B.C.  V2N 4Z9  CANADA
\endaddress

\email
walters\@hilbert.unbc.ca \ \ or \ \ walters\@unbc.ca  \hfill
\break\indent{\it Website}: http://hilbert.unbc.ca/walters
\endemail

\thanks
Research partly supported by NSERC grant OGP0169928
\hfill {\sevenrm (\TeX File: orthogonal.tex)}
\endthanks



\abstract 
Projections are constructed in the rotation algebra that are orthogonal 
to their Fourier transform and which are fixed under the flip automorphism.
Such projections are expected in a construction of an inductive limit structure
for the irrational rotation algebra that is invariant under the Fourier 
transform.  (Namely, as two circle algebras of the same dimension, that are
swapped by the Fourier transform, plus a few points.)
The calculation is based on Rieffel's construction of the Schwartz space as an 
equivalence bimodule of rotation algebras as well as on the theory of Theta functions. 
\endabstract

\endtopmatter

%%%%%%%%%%%%%%%%%%%%%%%%%%%%%  DOCUMENT BEGINS  %%%%%%%%%%%%%%%%%%%%%%%

\document

\specialhead \S1. Introduction \endspecialhead
Let $A_\theta$ denote the rotation C*-algebra (where $0<\theta<1$).  It is generated by
unitaries $U,V$ satisfying $VU=\lambda UV$ where $\lambda = e^{2\pi i\theta}$.
The Fourier transform is the automorphism $\sigma$ defined by
$\sigma(U)=V,\ \sigma(V)=U^{-1}$.  (Its square is the usual flip
automorphism.)  A hitherto open problem of Elliott, is whether 
$A_\theta$ is the inductive limit of Fourier invariant building blocks consisting 
of two circle algebras and a few points, and further, whether the fixed point subalgebra
$A_\theta^\sigma$ is an AF-algebra.  In [\SWc], a clue is given by constructing a 
model of an inductive limit automorphism of order four that agrees with $\sigma$ on 
$K_1$ and such that the fixed point subalgebra is approximately finite dimensional.  
The model predicts the existence of projections that are orthogonal to their
Fourier transform and fixed under the flip.  
The study of such projections would seem crucial in attempting to
unravel the structure of the Fourier transform.  Therefore, in this paper we show how 
one can explicitly construct such projections for $0<\theta<0.2345$, rational or 
irrational.  (In more recent work [\SWd], the results of this paper are used in an 
attempt to directly study Elliott's problem.)  
The construction here is based on Rieffel's realization [\MRb] of the Schwartz
space as an equivalence bimodule of rotation algebras associated with lattices $D$ in 
locally compact Abelian groups.  
The projection to be constructed has the C*-algebra-valued inner product form 
$\Dinner hh$ for an appropriate Schwartz function $h$ that involves Theta functions. 
Our main result is thus the following.

\proclaim{Theorem 1.1}
Let $0<\theta<0.2345$.  There exists a projection $e$ of trace $\theta$ in the 
rotation algebra $A_\theta$ such that
$$
e\sigma(e)=0, \qquad \sigma^2(e)=e.
$$
In particular, $e+\sigma(e)$ is a Fourier invariant projection of trace $2\theta$.
\endproclaim

This has the following immediate consequence.

\proclaim{Corollary 1.2}
Let $0<\theta<1$.  Each number in the set
$$
\{(m^2+n^2)\theta+k: m,n,k\in\Bbb Z\} \cap (0,0.2345)
$$
is the trace of a projection $e$ in $A_\theta$ satisfying $e\sigma(e)=0$ and 
$\sigma^2(e)=e$.  Further, $e+\sigma(e)$ is a Fourier invariant projection. 
\endproclaim
\demo{Proof} 
Let $\alpha = (m^2+n^2)\theta + k$ for some integer $k$ be in the above set.
Let $U_1 = \lambda^{mn/2}U^nV^m$ and $V_1=\sigma(U_1)=\lambda^{mn/2}V^nU^{-m}$. Then
$$
V_1U_1 = \lambda^{mn} V^n U^{-m} U^n V^m = \lambda^{mn/2} \lambda^{n(n-m)}
U^{n-m}V^{n+m}
= \lambda^{m^2+n^2} U_1V_1,
$$
and it is clear that $\sigma(V_1)=U_1^{-1}$.  Therefore $\sigma$ leaves invariant the
rotation subalgebra $A_{\alpha}$ generated by $U_1,V_1$, and induces the Fourier 
transform on it.  By the Theorem there is a projection $e$ in $A_{\alpha}$ such that
$e\sigma(e)=0,\ \sigma^2(e)=e$, and $\vartau(e) = \alpha$.
From this one gets the $\sigma$-invariant projection $e+\sigma(e)$ of
twice the trace.  \qed
\enddemo

We recall that for any $\theta$ in $(0,1)$ Boca [\FB] (Corollary 1.5) constructed 
a Fourier invariant projection in $A_\theta$ of trace $\theta$.  Presumably, the 
approach presented here combined with those of Boca (esp.~Section 2 of [\FB]) can 
be used to show that if $0 < \theta -\tfrac pq < \tfrac C{q^2}$, where $C$ is a 
positive constant akin to Boca's 0.948 (see his Proposition 2.1), or perhaps a 
little smaller, and if $p$ is a quadratic residue mod $q$, then there exists a 
projection of trace $q\theta-p$ satisfying the properties of Theorem 1.1 above.
(Similarly, for rational pairs satisfying $0 < \tfrac pq - \theta < \tfrac C{q^2}$ 
where $q-p$ is a quadratic residue mod $q$, one gets such a projection of trace 
$p-q\theta$.)

Throughout, we adopt the notation $e(t) := e^{2\pi it}$.

The technical part of Theorem 1.1 follows from the following more general result.

\proclaim{Theorem 1.3}
Let $A$ be a unital C*-algebra and $u,v$ be any pair of unitaries in $A$, where the spectrum
of $u$ is the unit circle.  Let $\alpha$ and $\beta$ be positive numbers satisfying 
$\beta^2 = 4(\alpha^2+1)$.  Then for $\alpha > 0.2568$ the following element
$$
X = \sum_{m,n\in\Bbb Z} e(\beta^2 mn/2) e^{-\pi \alpha\beta^2 n^2/2} 
e^{-\pi \beta^2 m^2/(2\alpha)} \Theta(m,n) \, u^m v^n
$$
is invertible in $A$, where
$$
\Theta(m,n) =
\vartheta_2(\tfrac\pi2 \beta^2 n, 2i\alpha)
\vartheta_3( i\tfrac{\pi}{2\alpha}\beta^2 m, i t_\alpha) +
\vartheta_3(\tfrac\pi2 \beta^2 n, 2i\alpha)
\vartheta_2( i\tfrac{\pi}{2\alpha}\beta^2 m, i t_\alpha)
$$
and $t_\alpha = 4\alpha + \tfrac2\alpha$.
\endproclaim

The Theta functions appearing here are recalled briefly in the next section.

\remark{Remark}
The commutation relation between $u,v$ is not needed to obtain the invertibility.  But it
is not clear if in general $X$ is positive, which it in fact is when $u,v$ are the
canonical unitary generators of the rotation algebra.  Theorem 1.3 holds if $v^n$ is
replaced by any sequence $v_n$ of unitaries in $A$ (while $u$ still has spectrum equal to
the unit circle).
\endremark

\bigpagebreak

%%%%%%%%%%%%%%%%%%%%%%%%%%%%%%%%%%%%%%%%%%%%%%%%%%%%%
%%%%%%%%%%%%%%%%%%%%%%%%%%%%%%%%%%%%%%%%%%%%%%%%%%%%%
\specialhead \S2.  Preliminaries \endspecialhead

Let us briefly recall Rieffel's setup in [\MRb].  Let $M$ be a locally compact Abelian
group and let $G=M\times \hat M$, where $\hat M$ is the dual of $M$.
Let $\frak h$ denote the Heisenberg cocycle on the group $G=M\times \hat M$ given by
$\frak h((m,s),(m',s')) = \langle m, s' \rangle$, where $\langle m, s' \rangle$ is the
canonical pairing $M\times\hat M \to \Bbb T$.  
The Heisenberg (unitary) representation $\pi:G\to \Cal U(L^2(M))$ of $G$ is given by
$$
(\pi_{(m,s)} f)(n) \ = \ \langle n,s \rangle f(n+m)
$$
where $m,n\in M,\  s\in\hat M$, and $f\in L^2(M)$.  It has the properties
$$
\pi_x\pi_y = \frak h(x,y) \pi_{x+y} = \frak h(x,y) \conj{\frak h(y,x)} \pi_y\pi_x, \qquad
\pi_x^* = \frak h(x,x) \pi_{-x}
$$
for $x,y\in G$.  If $D$ is a given lattice in $G$ (discrete cocompact subgroup), 
its covolume $|G/D|$ is the Haar measure of a fundamental domain for $D$.
Its associated twisted group C*-algebra $C^*(D,\frak h)$ is the 
C*-subalgebra of the bounded operators on $L^2(M)$ generated by the unitaries 
$\pi_x$ for $x\in D$.  The complementary lattice of $D$ is
$$
D^\perp = \{ y\in G: \frak h(x,y) \conj{\frak h(y,x)} = 1,\ \forall x\in D\}.
$$
The C*-algebra $C^*(D^\perp,\bar\frak h)$ can be viewed as the 
C*-subalgebra of bounded operators on $L^2(M)$ generated by the unitaries 
$\pi_y^*$ for $y\in D^\perp$.
The Schwartz space on $M$, denoted $\Cal S(M)$, is an equivalence bimodule with 
$C^*(D,\frak h)$ acting on the left and $C^*(D^\perp,\bar\frak h)$ acting on the 
right by 
$$
\align
af &= \int_D a(x)\pi_x(f) dx = |G/D|\,\sum_{x\in D} a(x)\pi_x(f)
\\
f b &= \int_{D^\perp} b(y) \pi_y^*(f) dy = \sum_{y\in D^\perp} b(y) \pi_y^*(f)
\endalign
$$
where $f\in\Cal S(M),\ a\in C^*(D,\frak h),\ b\in C^*(D^\perp,\bar\frak h)$, 
and where the mass point measure ($dx$) on $D$ is $|G/D|$ and on $D^\perp$ it is one.
The inner products on $\Cal S(M)$ with values in the algebras $C^*(D,\frak h)$ 
and $C^*(D^\perp, \conj{\frak h})$ are given, respectively, by
$$
\Dinner f g = |G/D|\,  \sum_{w\in D} \Dinner f g (w) \, \pi_w,
\qquad
\Dpinner f g = \sum_{z\in D^\perp} \Dpinner f g (z) \, \pi_z^*
$$
where
$$
\align
\Dinner f g (w_1,w_2) &= \int_M f(x) \conj{g(x+w_1)}\ \conj{\inner x {w_2}{}} dx 
\\
\Dpinner f g (z_1,z_2) &= \int_M \conj{f(x)} g(x+z_1) \inner x {z_2}{} dx
\endalign
$$
where $(w_1,w_2) \in D$ and $(z_1,z_2)\in D^\perp$.  These satisfy the associativity 
condition
$$
\Dinner f g h = f \Dpinner g h.
$$
(See [\MRb], pages 266 and 269.) By Rieffel's Theorem 2.15 in [\MRb], the Schwartz 
space $\Cal S(M)$ is an equivalence $C^*(D,\frak h)$-$C^*(D^\perp,\bar\frak h)$ bimodule.
The canonical normalized traces are given by
$$
\vartau_D\left( \sum_{w\in D} a_w\pi_w \right) = a_0,
\qquad
\vartau_{D^\perp}\left( \sum_{z\in D^\perp} b_z\pi_z^* \right) = b_0,
$$
($a_w,b_w\in \Bbb C$) which satisfy the equation
$$
\vartau_D( \Dinner f g ) = |G/D|\, \vartau_{D^\perp} (\Dpinner g f).
$$
From this it follows that if $f$ is a Schwartz function such that $\Dpinner f f = 1$,
then $e=\Dinner f f$ is a projection in $C^*(D,\frak h)$ of trace $|G/D|$.  
In this case, one has the isomorphism (which may be called the ``Morita transform")
$$
\mu: C^*(D^\perp,\bar\frak h) \to e C^*(D,\frak h) e, \qquad
\mu(x) = \Dinner {fx} f, \quad\text{and}\quad
\mu^{-1}(pyp) = \Dpinner f {yf}
$$
where $y\in C^*(D,\frak h)$ and $x\in C^*(D^\perp,\bar\frak h)$.
Projections arising in this manner we shall call {\it generalized Rieffel projections}.
This seems appropriate since in [\EL], Elliott and Lin showed that the classical
Rieffel projections [\MRa] in the rotation algbera can in fact be written in this 
manner.  (And perhaps also because of Proposition 2.8 of [\MRa].)

In the present paper, we shall take (as in [\FB]) $M=\Bbb R$ so that $\hat M
= M$ in a natural fashion.  This permits one to define the order four
automorphism $R:G\to G$ by $R(u;v)=(-v;u),\ u,v\in M$.  If $D$ is a lattice
subgroup of $G$ such that $R(D)=D$ (and hence $R(D^\perp)=D^\perp$), then (as
in [\FB]) there are order four automorphisms $\sigma, \sigma'$ of 
$C^*(D,\frak h)$ and $C^*(D^\perp,\bar \frak h)$, respectively, that satisfy
$$
\sigma(\pi_w) = \conj{\frak h(w,w)} \pi_{Rw}, \qquad
\sigma'(\pi_z^*) = \frak h(z,z) \pi_{Rz}^*,
\tag2.1
$$
for $w\in D,\ z\in D^\perp$, and
$$
\sigma(\Dinner fg) = \Dinner{\hat f}{\hat g}, \qquad
\sigma'(\Dpinner fg) = \Dpinner{\hat f}{\hat g},
\tag2.2
$$
where $\hat f$ is the Fourier transform of $f\in \Cal S(M)$ given by
$$
\hat f (s) = \int_M\ f(x) \conj{\inner x s {}}\ dx = 
\int_{-\infty}^\infty f(x) e(-sx) dx
$$
for $s\in M$.

The main Theta functions used in this paper are 
$$
\vartheta_2(z,t) = \sum_n \e^{\pi it(n+\frac 1 2)^2} \e^{i2z(n+\tfrac12)}, \qquad
\vartheta_3(z,t) = \sum_n \e^{\pi itn^2} \e^{i2zn}
$$
for $z,t\in\Bbb C$ and $\roman{Im}(t)>0$, where all summations range over the integers.
The zeros of $\vartheta_3$ are well-known to be given by
$$
\left( \frac \pi 2 + m\pi + \left( \frac \pi 2 + n\pi \right)t,\ t \right),
$$
where $m,n\in\Bbb Z$ and $\roman{Im}(t)>0$ are arbitrary.
We shall make use of the following identities for any integer $k$:
$$
\align
\vartheta_3(w + \tfrac\pi2 k it, it)
&= e^{-ikw} e^{\pi tk^2/4}\ \vartheta_{3-\bar k}(w, it)
\tag2.3
\\
\vartheta_2(w + \tfrac\pi2 k it, it)
&= e^{-ikw} e^{\pi tk^2/4}\ \vartheta_{2+\bar k}(w, it)
\tag2.4
\\
\vartheta(z+\pi it k, it) &= e^{-2kiz}e^{\pi t k^2} \vartheta(z, it)
\tag2.5
\endalign
$$
where $\bar k=0$ if $k$ is even and $\bar k = 1$ if $k$ is odd, and the third of these
holds for $\vartheta=\vartheta_2,\vartheta_3$.

Finally, we shall make use of the following equalities
$$
\int_{\Bbb R}\ e(Ax)\,e^{-\pi\alpha x^2}\,dx \ = \
\frac1{\sqrt\alpha} e^{-\pi A^2/\alpha}
\tag2.6
$$
and
$$
\iint_{\Bbb R^2}\ e(Ax+By)\,e^{-\pi\alpha(x^2+y^2)}\,e(-xy)\,dxdy \ = \
\tfrac1{\sqrt{\alpha^2+1}}
\exp\left(-\tfrac\pi{\alpha^2+1} [\alpha A^2+\alpha B^2 -2iAB] \right)
\tag2.7
$$
where $\alpha>0$ and $A,B\in\Bbb C$.

\bigpagebreak

%%%%%%%%%%%%%%%%%%%%%%%%%%%%%%%%%%%%%%%%%%%%%%%%%%%%%%%%%%%%%%%%%%%%%%%%%
%%%%%%%%%%%%  
\specialhead \S3. Schwartz Theta functions \endspecialhead

Fix $0<\theta<\tfrac12$, take $M=\Bbb R,\ G=\Bbb R \times \Bbb R$, and consider the 
lattice
$$
D:\ \ \bmatrix \varepsilon_1 \\ \varepsilon_2 \endbmatrix
=
\bmatrix
\sqrt\theta & 0 \\
0 & \sqrt\theta
\endbmatrix.
$$
It's covolume in $G$ is $|G/D|=\theta$ and the associated C*-algebra 
$C^*(D,\frak h)=A_\theta$ is generated by $U_j=\pi_{\varepsilon_j}$ 
which satisfy $U_1U_2=e(\theta)U_2U_1$.  Letting $\beta=1/\sqrt\theta$, 
the complementary lattice is
$$
D^\perp:\ \ \bmatrix \delta_1 \\ \delta_2 \endbmatrix
=
\bmatrix
0 & \beta \\
\beta & 0 
\endbmatrix.
$$
Letting $V_j=\pi_{\delta_j}^*=\pi_{-\delta_j}$ one has 
$V_1V_2=e(-\tfrac1\theta)V_2V_1 = e(-\beta^2)V_2V_1$,
and these unitaries generate the C*-algebra $C^*(D^\perp,\bar\frak h) = A_{1/\theta}$.
The automorphism $\sigma$ of $A_\theta$ as in (2.1) is given by 
$\sigma(U_1)=U_2,\ \sigma(U_2)=U_1^{-1}$, which is the Fourier transform.
The primary task of this paper is to prove the following.

\bigpagebreak

\proclaim{Theorem 3.1} For $0 < \theta < 0.2345$, there exists an even Schwartz 
function $h$ on $\Bbb R$ such that 
\itemitem{(i)} $\Dpinner h {\hat h}=0$,
\itemitem{(ii)} $\Dpinner h h$ is invertible,
\itemitem{(iii)} $\Dpinner h h$ is fixed under the flip $\sigma^2$.

Consequently, the element $e=\Dinner {ha}{ha}$ is a flip-invariant smooth projection in 
$A_\theta$ of trace $\theta$ such that $e\sigma(e)=0$, where $a = (\Dpinner h h)^{-1/2}$. 
In particular, $e+\sigma(e)$ is a Fourier invariant projection of trace $2\theta$.
\endproclaim

(The latter follows from $\Dpinner {ha}{ha}=1$.)

\bigpagebreak

\remark{Remark}
As commented below (see the remark before Proposition 4.4), one can in fact show, with a 
little extra effort, that for $0 < \theta \le 0.2427$ the element $\Dpinner h h$ is 
invertible for the Schwartz function $h$ constructed here.  
However, for $ \theta > 0.2451$ the methods here would have to be modified.
\endremark

\bigpagebreak

We shall now examine the Schwartz functions on $\Bbb R$ to be used in the 
construction.  These functions have the Gaussian-Theta form
$$
g_{r, \gamma} (x) 
\ = \ 
e^{-\pi\alpha x^2}\ \sum_p\ a_p e([rp-\gamma]x) 
\ = \ 
e^{-\pi\alpha x^2}\ e(-\gamma x) \, \vartheta_3(-\tfrac12\pi i\alpha + \pi rx, i\alpha)
$$
where
$$
a_p = e^{-\pi \alpha p^2} e^{\pi\alpha p}
$$
and $\alpha>0$ and $r, \gamma$ are real.  It is not hard to check that 
$g = g_{r,\gamma}\in \Cal S(\Bbb R)$.  We will choose these parameters so that
$$
\Dpinner g {\hat g} = 0 \quad \text{and} \quad
\Dpinner {\tilde g} {\hat g} = 0,
\tag3.1
$$
where $\tilde g(x) = g(-x) = \hat{\hat g}(x)$.
One then puts $h=g+\tilde g$ so that $h$ is an even function satisfying 
$\Dpinner h {\hat h} = 0$ and $\Dinner h h$ is fixed under the flip automorphism.
For certain of these parameters it will be shown in the next section that 
$\Dpinner h h$ is invertible.

For real $r,s,\gamma,\gamma'$, one has
$$
\align
&\Dpinner {g_{r,\gamma}}{\hat g_{s,\gamma'}}(m\delta_1+n\delta_2) 
= \int_{\Bbb R}\ \conj{g_{r,\gamma}(x)} \hat g_{s,\gamma'}(x+\beta n) e(\beta mx)\, dx \\
&= \iint_{\Bbb R^2}\ \conj{g_{r,\gamma}(x)} g_{s,\gamma'}(y) e(-xy) e(\beta mx - \beta ny)\,dxdy\\
&= \sum_{p,q} \bar a_p a_q \ \iint_{\Bbb R^2}\ 
e([(\beta m - (rp-\gamma))x + ((sq-\gamma') - \beta n)y]\, e^{-\pi\alpha(x^2+y^2)}
e(-xy) \, dxdy \\
&= \tfrac1{\sqrt{\alpha^2+1}} \sum_{p,q} \bar a_p a_q
\exp\left( -\tfrac\pi{\alpha^2+1} D_{pq} \right)
\endalign
$$
using (2.7), where
$$
\align
D_{pq} &= \alpha [\beta m - rp + \gamma]^2 + \alpha[sq - \gamma' - \beta n]^2 - 
2i[\beta m - rp + \gamma]\,[sq - \gamma' - \beta n]
\\
& = \alpha(\beta m - rp)^2  + 2\alpha(\beta m - rp) \gamma + \alpha \gamma^2 
 + \alpha(sq - \beta n)^2  - 2\alpha(sq - \beta n) \gamma' + \alpha (\gamma')^2  \\
&\ \ \  - 2i(\beta m - rp)(sq-\beta n) + 2i(\beta m - rp)\gamma' - 2i \gamma(sq-\beta n) + 
2i\gamma\gamma' 
\\
& = \alpha r^2 p^2 - 2r(\alpha\beta m + i \beta n + \alpha\gamma + i \gamma') p 
+ \alpha s^2 q^2 - 2s(\alpha\beta n + i\beta m + \alpha \gamma' + i \gamma) q \\
& \ \ \ + 2irspq + K',
\endalign
$$
where $K'$ consists of all terms not involving $p,q$.  In order to make the cross term
involving $2irspq$ disappear in the summation, we shall assume that $|rs|=\alpha^2+1$.
This makes the above summation separable and one gets
$$
\Dpinner {g_{r,\gamma}}{\hat g_{s,\gamma'}}(m\delta_1+n\delta_2) = C_{mn} E_1(m,n) E_2(m,n)
$$
where $C_{mn}$ is an exponential constant (independent of $p,q$ but depending on $m,n$),
and where
$$
\align
E_1(m,n) &= \sum_p \bar a_p \exp\left( -\tfrac\pi{\alpha^2+1} [\alpha r^2 p^2 - 2r(\alpha\beta m
+ i \beta n + \alpha\gamma + i \gamma') p] \right),
\\
E_2(m,n) &= \sum_q a_q \exp\left( -\tfrac\pi{\alpha^2+1} [\alpha s^2 q^2 - 2s(\alpha\beta n +
i\beta m + \alpha \gamma' + i \gamma) q] \right).
\endalign
$$
Substituting $a_p$, one has
$$
\align
E_1(m,n) &= \sum_p e^{-\pi d p^2} 
\exp\left( \pi\alpha p + \tfrac{2\pi r}{\alpha^2+1} [\alpha\beta m + i \beta n +
\alpha\gamma + i \gamma']\, p \right)
\\
&= \vartheta_3\left( \tfrac12 i\pi \alpha + \tfrac{i\pi r}{\alpha^2+1} [\alpha\beta m + i
\beta n + \alpha\gamma + i \gamma'],\ id \right)
\endalign
$$
where $d = \alpha + \tfrac{\alpha r^2}{\alpha^2+1}$.  Similarly, one can write $E_2$.
In order to satisfy (3.1), one can choose either to make $E_1\equiv0$ or 
$E_2\equiv0$.  For our purposes, it will be enough to make $E_1\equiv0$ for 
$\gamma'=\pm\gamma$.  This will be arranged as follows.  
Since $\tilde g_{r,\gamma} = g_{-r,-\gamma}$, we want to have
$$
\Dpinner{g_{r, \gamma}}{\hat g_{r,\gamma}} = 0, \qquad \text{and} \qquad
\Dpinner{g_{r, \gamma}}{\hat g_{-r,-\gamma}} = 0,
$$
and in both cases we have $s=\pm r,\ \gamma'=\pm\gamma$.  Thus $r^2=s^2=\alpha^2+1$, 
and $d=2\alpha$.  
Now $E_1\equiv0$ iff for each pair of integers $m,n$ there are integers $M,N$ such that
$$
\frac12 i \alpha + \frac{ir}{\alpha^2+1} [\alpha\beta m + i
\beta n + \alpha\gamma + i \gamma'] 
= \frac12 + M + \left( \frac12 + N \right) i 2\alpha.
$$
Putting $m=n=0$ gives integers $M_0,N_0$ such that
$$
\frac12 i \alpha + \frac{ir}{\alpha^2+1} [\alpha\gamma + i \gamma'] 
= \frac12 + M_0 + \left( \frac12 + N_0 \right) i 2\alpha.
\tag3.2
$$
and upon subtracting one gets
$$
\tfrac{r}{\alpha^2+1} [i\alpha\beta m - \beta n ]  =  M' + N'i2\alpha,
$$
where $M',N'$ are integers depending on $m,n$.  Thus
$$
- \tfrac{r\beta}{\alpha^2+1} n  =  M',
\qquad
\tfrac{r\beta}{\alpha^2+1} m  =  2N'.
$$
The first of these implies that $\tfrac{|r|\beta}{\alpha^2+1} = v$ must be a 
positive integer.
And from the second, one must have $v=2L\ge2$ for some positive integer $L$.  
Further, $r^2=s^2=\beta^2/v^2=\alpha^2+1$.
To get the greater range for $\theta = 1/\beta^2 \le 1/v^2$, we shall henceforth choose
$v=2$ ($L=1$).  Therefore, we are lead to make the following choices
$$
|r| = |s| = \frac \beta 2, \qquad \beta^2 = 4(\alpha^2+1).
$$
So we take $g := g_{\beta/2, \gamma}$.  Now to satisfy 
$\Dpinner{g_{\beta/2, \gamma}}{\hat g_{\beta/2, \gamma}} = 0$, 
condition (3.2) must be satisfied with $\gamma'=\gamma$ and becomes
$$
\frac12 i \alpha + \frac{\beta\gamma}{2(\alpha^2+1)} [i\alpha-1]
= \frac12 + M_0 + \left( \frac12 + N_0 \right) 2i\alpha
\tag3.3
$$
and in order to satisfy 
$\Dpinner{g_{\beta/2, \gamma}}{\hat g_{-\beta/2, -\gamma}} = 0$, 
condition (3.2) must be satisfied with $\gamma'=-\gamma$ and becomes
$$
\frac12 i \alpha + \frac{\beta\gamma}{2(\alpha^2+1)} [i\alpha+1]
= \frac12 + M_0' + \left( \frac12 + N_0' \right) 2i\alpha.
\tag3.4
$$
Subtracting (3.3) from (3.4) implies that $\tfrac{\beta\gamma}{\alpha^2+1} =
\tfrac{4\gamma}\beta$ must be an integer.  We choose it to be 1.  Hence we take
$$
\gamma = \frac \beta 4.
$$
(This choice will in fact turn out to facilitate our computations with Theta functions 
below and guarantees that $\Dpinner h h$ is invertible for $\beta$ a little greater 
than 2.) 
Now with $r=\tfrac\beta2$ and $\gamma=\tfrac\beta4$, (3.4) follows from (3.3) and the 
latter clearly holds with $M_0=-1,\,N_0=0$.

We therefore conclude that
$$
\Dpinner{g_{\beta/2,\beta/4}}{\hat g_{\beta/2,\beta/4}} = 0, \qquad
\Dpinner{g_{\beta/2,\beta/4}}{\hat g_{-\beta/2,-\beta/4}} = 0,
$$
as desired.  For simplicity, we shall write $g := g_{\beta/2,\beta/4}$ (where $\beta$
is fixed).  More specifically,
$$
g(x) = e^{-\pi\alpha x^2}\ \sum_p\ a_p e([\tfrac\beta2 p - \tfrac\beta4]x), \qquad
a_p = e^{-\pi\alpha p^2} e^{\pi\alpha p}.
$$

\remark{Remark} 
If we drop the $\tfrac\beta4$ term a calculation similar to that done in the next section
shows that $\Dpinner h h$ is singular.
\endremark

%%%%%%%%%%%%%%%%%%%%%%%%%%%%%%%%%%%%%%%%%%%%%%%%%%%%%%
%%%%%%%%%%%%%%%%%%%%%%%%%%%%%%%%%%%%%%%%%%%%%%%%%%%%%%
\vfill\eject

\specialhead \S4. Invertibility of $\Dpinner h h$ \endspecialhead
This section is aimed at showing that the positive element $\Dpinner h h$ is invertible,
where $h$ was constructed in the previous section.  The calculation here also proves 
Theorem 1.3 (including the last assertion in the remark following it).  
We begin with a lemma. 

\bigpagebreak

\proclaim{Lemma 4.1}
Let there be given two Gaussian-Theta Schwartz functions of the form 
$$
f_1(x) = e^{-\pi\alpha_1 x^2}\ \sum_p\ a_p e([r_1p-\gamma_1]x), \qquad
f_2(x) = e^{-\pi\alpha_2 x^2}\ \sum_q\ b_q e([r_2q-\gamma_2]x), 
$$
where $\alpha_j>0,\ r_j,\gamma_j$ are real, and $a_p,b_q$ are rapidly
decreasing sequences, and let $\alpha = \alpha_1+\alpha_2$.  For any real $s,t$, one has
$$
\multline
\int_{\Bbb R} \conj{f_1(x)}{f_2(x+s)} e(tx)\,dx
\\ =
\frac{e(-\alpha_2st/\alpha)}{\sqrt{\alpha}}
e^{-\pi \alpha_1\alpha_2 s^2/\alpha}
\ e^{-\pi t^2/\alpha}
e^{-\pi(\gamma_1-\gamma_2)^2/\alpha}
\cdot e^{2\pi(\gamma_2-\gamma_1)t/\alpha}
e(-(\alpha_2\gamma_1+\alpha_1\gamma_2)\,\tfrac s\alpha)
\Theta(s,t)
\endmultline
$$
where
$$
\multline
\Theta(s,t) = 
\sum_{p,q} \bar a_p b_q \ 
e([\alpha_1r_2q + \alpha_2r_1p]\,\tfrac s\alpha)
\exp\left(\tfrac{2\pi}\alpha(\gamma_2-\gamma_1)(r_2q-r_1p) \right)  \\
\cdot \exp\left(-\tfrac\pi\alpha [(r_2q-r_1p)^2 + 2t(r_2q-r_1p)]
\right)
\endmultline
$$
\endproclaim
\demo{Proof} We have
$$
\align
&\int_{\Bbb R}\ \conj{f_1(x)} f_2(x+s) e(tx)\, dx \\
&= 
\sum_{p,q} \bar a_p b_q \ \int_{\Bbb R}\
e^{-\pi\alpha_1x^2} e(-[r_1p-\gamma_1]x) e^{-\pi\alpha_2(x+s)^2} 
e([r_2q-\gamma_2](x+s)) e(tx)\,dx
\\
&=
\sum_{p,q} \bar a_p b_q e([r_2q-\gamma_2]s)\ 
\int_{\Bbb R}\ e^{-\pi[\alpha_1x^2+\alpha_2(x+s)^2]} 
e([r_2q-r_1p+\gamma_1-\gamma_2+t] x)\ dx
\\
&= 
\sum_{p,q} \bar a_p b_q \ e([r_2q-\gamma_2]s) 
e^{-\pi \alpha_1\alpha_2 s^2/\alpha}
\ \int_{\Bbb R}\
e([r_2q-r_1p+\gamma_1-\gamma_2+t]\,x) \ 
e^{-\pi\alpha\left(x+\tfrac{\alpha_2 s}\alpha \right)^2}\ dx
\\
&=
\sum_{p,q} \bar a_p b_q \ e([r_2q-\gamma_2]s)
e^{-\pi \alpha_1\alpha_2 s^2/\alpha}
e(-[r_2q-r_1p+\gamma_1-\gamma_2+t]\tfrac{\alpha_2 s}\alpha)
\\ 
&\hskip2in 
\cdot \int_{\Bbb R}\
e([r_2q-r_1p+\gamma_1-\gamma_2+t]\,x) \ e^{-\pi\alpha x^2}\ dx
\\
&=
\frac1{\sqrt\alpha} e^{-\pi \alpha_1\alpha_2 s^2/\alpha}
e(-[\gamma_1-\gamma_2+t]\tfrac{\alpha_2s}\alpha)
\\ 
&\hskip0.4in \cdot 
\sum_{p,q} \bar a_p b_q \ e([r_2q-\gamma_2]s)
e(-[r_2q-r_1p]\tfrac{\alpha_2s}\alpha)
\exp\left(-\tfrac\pi\alpha [r_2q-r_1p+\gamma_1-\gamma_2+t]^2 \right) 
\endalign
$$
and since
$$
\multline
[r_2q-r_1p+\gamma_1-\gamma_2+t]^2 \ = \ 
(r_2q-r_1p)^2 + 2t (r_2q-r_1p) + t^2 \\
+ 2(\gamma_1-\gamma_2)(r_2q-r_1p) + 2t(\gamma_1-\gamma_2) + 
(\gamma_1-\gamma_2)^2
\endmultline
$$
one gets
$$
\align
\int_{\Bbb R}\ \conj{f_1(x)} & f_2(x+s) e(tx)\, dx
\\
&=
\frac1{\sqrt\alpha} e^{-\pi \alpha_1\alpha_2 s^2/\alpha}
e(-[\gamma_1-\gamma_2+t]\tfrac{\alpha_2s}\alpha)
\exp\left( -\tfrac\pi\alpha(\gamma_1-\gamma_2)^2 \right) 
\exp\left( 2\tfrac\pi\alpha t(\gamma_2-\gamma_1) \right)
\\
&\ \ \  \cdot
e^{-\tfrac\pi\alpha t^2}
\sum_{p,q} \bar a_p b_q \ e([r_2q-\gamma_2]s)
e(-[r_2q-r_1p]\tfrac{\alpha_2s}\alpha)
\exp(\tfrac{2\pi}\alpha(\gamma_2-\gamma_1)[r_2q-r_1p])
\\
&\hskip2in \cdot 
\exp\left( -\tfrac\pi\alpha[(r_2q-r_1p)^2 + 2t(r_2q-r_1p) \right)
\endalign
$$
which gives the result in the statement of the lemma. \qed
\enddemo

\bigpagebreak

%%%%%%%%%%%%%%%%%%%%%%%%%%%%%%%%%%%%%%%%%%%%%%%%%%%%%
%%%%%%%%%%%%%%%%%%%%%%%%%%%%%%%%%%%%%%%%%%%%%%%%%%%%%

%%\newpage

Applying Lemma 4.1 to the functions $g_\rho := g_{\rho\beta/2, \rho\beta/4}$ 
for $\rho=\pm1$, one obtains (where $\rho,\nu=\pm1$)
$$
\multline
\int_{\Bbb R}\ \conj{g_\rho(x)} g_\nu(x+s) e(tx)\, dx
\\
= \frac{e(-\tfrac{st}2)}{\sqrt{2\alpha}}
e^{-\pi \alpha s^2/2} e^{-\pi t^2/(2\alpha)}
\cdot
\exp\left( -\tfrac{\pi\beta^2}{32\alpha}(\rho-\nu)^2 \right)
\exp\left( \tfrac{\pi\beta t}{4\alpha}(\nu-\rho) \right)
e(-\tfrac{\beta s}8(\rho+\nu))\ \Theta_{\rho,\nu}
\endmultline
$$
where 
$$
\multline
\Theta_{\rho,\nu} = 
\sum_{p,q} \ e^{-\pi\alpha(p^2+q^2)}\ e^{\pi\alpha(p+q)}\ 
e(\tfrac{\beta}4 [\rho p + \nu q]\,s)
\exp\left( -\tfrac{\pi\beta^2}{8\alpha}(\rho-\nu)(\nu q - \rho p) \right)
\\ 
\exp\left(-\tfrac{\pi\beta}{8\alpha} [\beta(\nu q - \rho p)^2 + 4t(\nu q - \rho p)] 
\right).
\endmultline
$$
Making the replacements $p\to \rho p,\ q\to \nu q$ this becomes 
$$
\multline
\Theta_{\rho,\nu} =
\sum_{p,q} \ e^{-\pi\alpha(p^2+q^2)}\ e^{\pi\alpha(\rho p + \nu q)}\
e(\tfrac{\beta}4 (p+q)\, s)
\exp\left( -\tfrac{\pi\beta^2}{8\alpha}(\rho-\nu)(q-p) \right)
\\
\exp\left(-\tfrac{\pi\beta}{8\alpha} [\beta(q-p)^2 + 4t(q-p)] \right).
\endmultline
$$
Now let $k=q-p$:
$$
\multline
\Theta_{\rho,\nu} =
\sum_{p,k} \  e^{-\pi\alpha[p^2+(p+k)^2]} e^{\pi\alpha[\rho p + \nu(p+k)]} 
e\left( \tfrac{\beta}4(k+2p) s \right)\
\exp\left( -\tfrac{\pi\beta^2}{8\alpha}(\rho-\nu)\,k \right)
\\
\exp\left(-\tfrac{\pi\beta}{8\alpha} [\beta k^2 + 4tk] \right).
\endmultline
$$
To calculate this more explicitly, we write
$$
\Theta_{\rho,\nu} = \sum_k e^{-\pi t_\alpha k^2/4}
\, \exp\left( [\pi\alpha\nu - \tfrac{\pi\beta t}{2\alpha} 
	            - \tfrac{\pi\beta^2(\rho-\nu)}{8\alpha} ]\,k \right) \ 
e^{-\pi\alpha (\rho+\nu) k/2} \ H(k)
$$
where $t_\alpha = 4\alpha + \tfrac2\alpha = \tfrac{\beta^2}{2\alpha} + 2\alpha$ and
$$
\align
H(k) &= \sum_p e^{-2\pi\alpha(p+\tfrac k2)^2}\, 
\exp\left( [\pi i\beta s + \pi\alpha(\rho+\nu)] (p+\tfrac k2) \right)
\\
&= \sum_p e^{-2\pi\alpha(p+\tfrac k2)^2}\, e^{L (p+\tfrac k2)},
\endalign
$$
where
$$
L = L_{\rho,\nu} := \pi i\beta s + \pi\alpha(\rho+\nu).
$$
For $k$ even, one has $H(k) = \vartheta_3(\tfrac i2L, 2i\alpha)$, and for $k$ odd 
$H(k) = \vartheta_2(\tfrac i2L, 2i\alpha)$.  Letting 
$$
M = M_{\rho,\nu} = \pi\alpha\nu - \tfrac{\pi\beta t}{2\alpha} - 
\tfrac{\pi\beta^2(\rho-\nu)}{8\alpha} - \tfrac\pi2\alpha(\rho+\nu)
\ = \ 
\tfrac\pi2\alpha(\nu-\rho) -
\tfrac{\pi\beta t}{2\alpha} + \tfrac{\pi\beta^2(\nu-\rho)}{8\alpha},
$$
which, from $\beta^2=4(\alpha^2+1)$, becomes
$$
M = \tfrac\pi4 t_\alpha (\nu-\rho) - \tfrac{\pi\beta t}{2\alpha},
$$
and one obtains
$$
\align
\Theta_{\rho,\nu} &= \sum_k e^{-\pi t_\alpha k^2/4} e^{Mk} H(k) 
\\
&= \vartheta_3(\tfrac i2L, 2i\alpha) \sum_k e^{-\pi t_\alpha k^2} e^{2Mk} +
\vartheta_2(\tfrac i2L, 2i\alpha) \sum_k e^{-\pi t_\alpha (k+\tfrac12)^2} e^{2M(k+\tfrac12)}
\\
&= 
\vartheta_3(\tfrac i2L, 2i\alpha) \vartheta_3(iM, it_\alpha) +
\vartheta_2(\tfrac i2L, 2i\alpha) \vartheta_2(iM, it_\alpha).
\endalign
$$
Let $\ell = \tfrac12(\nu+\rho)$ and $d = \tfrac12(\nu-\rho$) which are equal to $0,\pm1$.
We have thus obtained
$$
\align
\int_{\Bbb R}\ \conj{g_\rho(x)} g_\nu(x+s) e(tx)\, dx
&=
\frac{e(-\tfrac12 st)}{\sqrt{2\alpha}}
e^{-\pi \alpha s^2/2} e^{-\pi t^2/(2\alpha)}
e^{-\tfrac{\pi\beta^2}{8\alpha} d^2}
e^{\tfrac{\pi\beta t}{2\alpha}d}
e^{-i\pi\tfrac{\beta s}2 \ell}
\\
&\ \ \ \ \ \ \ \cdot
\Bigl[
\vartheta_3(\tfrac i2L, 2i\alpha) \vartheta_3(iM, it_\alpha) +
\vartheta_2(\tfrac i2L, 2i\alpha) \vartheta_2(iM, it_\alpha)
\Bigr]. 
\endalign
$$
Using the identities (2.3) and (2.4) each of these Theta functions can be simplified 
as follows
$$
\align
\vartheta_3(\tfrac i2L, 2i\alpha) 
&=
e^{\tfrac{i\pi}2\beta s\ell} e^{\tfrac\pi2 \alpha\ell^2}
\vartheta_{3-\bar\ell}(\tfrac\pi2 \beta s, 2i\alpha)
\\
\vartheta_2(\tfrac i2L, 2i\alpha) 
&=
e^{\tfrac{i\pi}2\beta s\ell} e^{\tfrac\pi2 \alpha\ell^2}
\vartheta_{2+\bar\ell}(\tfrac\pi2 \beta s, 2i\alpha)
\\
\vartheta_3(iM, it_\alpha)
&=
e^{-\tfrac{\pi}{2\alpha}d\beta t} e^{\tfrac\pi4 d^2 t_\alpha}
\vartheta_{3-\bar d}( i\tfrac{\pi}{2\alpha}\beta t, i t_\alpha)
\\
\vartheta_2(iM, it_\alpha)
&=
e^{-\tfrac{\pi}{2\alpha}d\beta t} e^{\tfrac\pi4 d^2 t_\alpha}
\vartheta_{2+\bar d}( i\tfrac{\pi}{2\alpha}\beta t, i t_\alpha)
\endalign
$$
and since $d^2+\ell^2 = \tfrac12(\nu^2+\rho^2) = 1$ we get
$$
\align
\int_{\Bbb R}\ \conj{g_\rho(x)} g_\nu(x+s) e(tx)\, dx
&=
\frac{e(-\tfrac12 st)}{\sqrt{2\alpha}}
e^{-\pi \alpha s^2/2} 
e^{-\pi t^2/(2\alpha)}
e^{\tfrac\pi2 \alpha}
\\
&\cdot
\Bigl[
\vartheta_{3-\bar\ell}(\tfrac\pi2 \beta s, 2i\alpha)
\vartheta_{3-\bar d}( i\tfrac{\pi}{2\alpha}\beta t, i t_\alpha)
+
\vartheta_{2+\bar\ell}(\tfrac\pi2 \beta s, 2i\alpha)
\vartheta_{2+\bar d}( i\tfrac{\pi}{2\alpha}\beta t, i t_\alpha)
\Bigr].
\endalign
$$
Since $h=g_1 + g_{-1}$, summing these over $\rho,\nu=\pm1$ yields the following.

%%%%%%%%%%%%%%%%%%%%%%%%%%%%%%%%%%%%%%%%%%%
\proclaim{Lemma 4.2} For real $s,t$ one has
$$
\int_{\Bbb R}\ \conj{h(x)} h(x+s) e(tx)\, dx
=
\frac{e(-\tfrac12 st)}{\sqrt{2\alpha}}
e^{-\pi \alpha s^2/2} 
e^{-\pi t^2/(2\alpha)}
\Gamma(\tfrac t\beta, \tfrac s\beta)
$$
where
$$
\Gamma(u,v) =
4e^{\tfrac\pi2 \alpha}
\Bigl[
\vartheta_2(\tfrac\pi2 \beta^2 v, 2i\alpha)
\vartheta_3( i\tfrac{\pi}{2\alpha}\beta^2 u, i t_\alpha)
+
\vartheta_3(\tfrac\pi2 \beta^2 v, 2i\alpha)
\vartheta_2( i\tfrac{\pi}{2\alpha}\beta^2 u, i t_\alpha)
\Bigr]
$$
for $u,v\in\Bbb R$.
\endproclaim
This lemma now gives
$$
\Dpinner h h (m\delta_1+n\delta_2) = 
\tfrac1{\sqrt{2\alpha}} e(-\beta^2 mn/2) e^{-\pi\alpha\beta^2 n^2/2} 
e^{-\pi\beta^2m^2/(2\alpha)} \Gamma(m,n)
$$
hence
$$
\align
\Dpinner h h &= \sum_{m,n} \Dpinner h h (m\delta_1+n\delta_2) V_2^n V_1^m 
\\
&=  \sum_{m,n} e(\beta^2 mn) \Dpinner h h (m\delta_1+n\delta_2) V_1^m V_2^n
\\
&= \tfrac1{\sqrt{2\alpha}}
\sum_{m,n} e(\beta^2 mn/2) e^{-\pi \alpha\beta^2 n^2/2} e^{-\pi \beta^2 m^2/(2\alpha)}
\Gamma(m,n) \, V_1^m V_2^n
\\
&= \tfrac1{\sqrt{2\alpha}}
\sum_n \ e^{-\pi \alpha\beta^2 n^2/2} \, \psi_n V_2^n
\endalign
$$
where $\psi_n = \psi_n(V_1)$, viewed as a period 1 function on $\Bbb R$, is given by
$$
\align
\psi_n(t) &:= \sum_m\ e(\beta^2 mn/2) e^{-\pi \beta^2 m^2/(2\alpha)}\Gamma(m,n)\, e(mt)
\\
&=
4e^{\tfrac\pi2 \alpha} 
\Bigl[ 
\vartheta_2(\tfrac\pi2 \beta^2 n, 2i\alpha) F_n(t') +
\vartheta_3(\tfrac\pi2 \beta^2 n, 2i\alpha) G_n(t')
\Bigr]
\endalign
$$
where $t' = t + \tfrac12\beta^2 n$ and 
$$
F_n(t) = \sum_m\ e^{-\pi \beta^2 m^2/(2\alpha)} 
\vartheta_3( i\tfrac{\pi}{2\alpha}\beta^2 m, i t_\alpha)\, e(mt)
$$
and
$$
G_n(t) = \sum_m\ e^{-\pi \beta^2 m^2/(2\alpha)} 
\vartheta_2( i\tfrac{\pi}{2\alpha}\beta^2 m, i t_\alpha)\, e(mt).
$$
Note that up to a multiplicative positive constant, $\Dpinner h h$ is the element $X$ of
Theorem 1.3, whose proof is furnished by the following argument.
We now calculate $F_n$ and $G_n$ explicitly as follows.  For convenience, let
$$
\tau_\alpha := \tfrac{\beta^2}{2\alpha} = 2(\alpha + \tfrac1\alpha) = 
t_\alpha - 2\alpha.
$$
Inserting the series of $\vartheta_3$ one has
$$
\align
F_n(t) &= \sum_m e^{-\pi \beta^2 m^2/(2\alpha)}
\sum_k e^{-\pi t_\alpha k^2} e^{2ik[i\tfrac{\pi}{2\alpha}\beta^2 m]} 
\,e^{2\pi imt}
\\
&= \sum_k e^{-\pi t_\alpha k^2}
\sum_m e^{-\pi \beta^2 m^2/(2\alpha)} e^{2im[\pi t + i\tfrac{\pi}{2\alpha}\beta^2 k]}
\\
&= \sum_k e^{-\pi t_\alpha k^2}
\sum_m e^{-\pi \tau_\alpha m^2} e^{2im[\pi t + i\pi\tau_\alpha k]}
\\
&= \sum_k e^{-\pi t_\alpha k^2}
e^{\pi\tau_\alpha k^2} \sum_m e^{-\pi \tau_\alpha (m+k)^2} e^{2im \pi t}
\\
&= 
\sum_k e^{-2\pi\alpha k^2} e^{-2ik \pi t} \sum_m e^{-\pi \tau_\alpha m^2} e^{2im \pi t} 
\\
&=
\vartheta_3(\pi t, 2i\alpha)\, \vartheta_3(\pi t, i\tau_\alpha)
\endalign
$$
A similar computation gives
$$
G_n(t) = \vartheta_2(\pi t, 2i\alpha)\, \vartheta_2(\pi t, i\tau_\alpha). 
$$
Thus we obtain
$$
\psi_n(t) = 4e^{\pi\alpha/2} \Bigl[
\vartheta_2(v_n, 2i\alpha)\vartheta_3(\pi t',2i\alpha)\vartheta_3(\pi t',i\tau_\alpha)
+
\vartheta_3(v_n, 2i\alpha)\vartheta_2(\pi t',2i\alpha)\vartheta_2(\pi t',i\tau_\alpha)
\Bigr].
$$
where $v_n := \tfrac\pi2 \beta^2 n$ and $t' := t + \tfrac12\beta^2n$.
Write this in the form
$$
\psi_n(t) = 4e^{\pi\alpha/2} 
\vartheta_3(v_n, 2i\alpha)\, \vartheta_3(\pi t', 2i\alpha)\,
\vartheta_3(\pi t', i\tau_\alpha) \cdot Q_n(t)
$$
where
$$
Q_n(t) = \frac{\vartheta_2(v_n, 2i\alpha)}{\vartheta_3(v_n, 2i\alpha)} + 
\frac{\vartheta_2(\pi t', 2i\alpha)}{\vartheta_3(\pi t', 2i\alpha)} 
\frac{\vartheta_2(\pi t', i\tau_\alpha)}{\vartheta_3(\pi t', i\tau_\alpha)}.
$$
It is clear that $\psi_n$ is a real function.  Our next goal is to show that $\psi_0$ is
positive.  We have
$$
\psi_0(t) =
4e^{\pi\alpha/2} \vartheta_3(0, 2i\alpha)\, \vartheta_3(\pi t, 2i\alpha)\,         
\vartheta_3(\pi t, i\tau_\alpha) \cdot Q_0(t),
$$
$$
Q_0(t) = \frac{\vartheta_2(0, 2i\alpha)}{\vartheta_3(0, 2i\alpha)} + 
\frac{\vartheta_2(\pi t, 2i\alpha)}{\vartheta_3(\pi t, 2i\alpha)} 
\frac{\vartheta_2(\pi t, i\tau_\alpha)}{\vartheta_3(\pi t, i\tau_\alpha)}
$$
and it is enough to show that $Q_0(t)>0$.  This is an immediate consequence of the 
following two facts:

\itemitem{(F1)} For all real $s,t$ and $\alpha>0$, one has 
$$
\frac{\vartheta_2(\pi t, i\tau_\alpha)}{\vartheta_3(\pi s, i\tau_\alpha)} \ < \ 0.09.
$$

\itemitem{(F2)} For each real $t$, one has
$$
\left| \frac{\vartheta_2(\pi t, 2i\alpha)}{\vartheta_3(\pi t, 2i\alpha)} \right| \ \le \
\frac{\vartheta_2(0, 2i\alpha)}{\vartheta_3(0, 2i\alpha)}.
$$

The first of these follows from the inequality 
$\vartheta_2(0,i\tau_\alpha) \le \vartheta_2(0,4i)$ since 
$\tau_\alpha = 2(\alpha + \tfrac1\alpha) \ge4$ for all 
$\alpha>0$, and from the estimate 
$\vartheta_3(\tfrac\pi2, i\tau_\alpha) > 1 - 2e^{-\pi\tau_\alpha} \ge 1 - 2e^{-4\pi}$.
Thus for real $s,t$,
$$
\frac{\vartheta_2(\pi t, i\tau_\alpha)}{\vartheta_3(\pi s, i\tau_\alpha)} \ \le
\frac{\vartheta_2(0, i\tau_\alpha)}{\vartheta_3(\tfrac\pi2, i\tau_\alpha)} \ 
< \ \frac{\vartheta_2(0,4i)}{1 - 2e^{-4\pi}} \approx 0.0864
$$
which proves (F1).  (Note that here we have used the well-known fact that for each real 
$t$ and $x>0$ one has
$0 < \vartheta_3(\tfrac\pi2, ix) \le \vartheta_3(t, ix) \le \vartheta_3(0,ix)$.)
To show (F2), one need only check that the function
$t \mapsto \tfrac{\vartheta_2(\pi t, 2i\alpha)}{\vartheta_3(\pi t, 2i\alpha)}$
has peroid $2$, is decreasing over the interval $[0,1]$, increasing over $[1,2]$, and 
hence attains its maximum value when $t=0$.

\bigpagebreak

We have
$$
\left| \frac{\psi_n(t)}{\psi_0(s)} \right|
= 
\frac{\vartheta_3(v_n, 2i\alpha)}{\vartheta_3(0, 2i\alpha)}\, 
\frac{\vartheta_3(\pi t', 2i\alpha)}{\vartheta_3(\pi s, 2i\alpha)}\,
\frac{\vartheta_3(\pi t', i\tau_\alpha)}{\vartheta_3(\pi s, i\tau_\alpha)} 
\cdot 
\frac{Q_n(t)}{Q_0(s)}.
$$
We make an upper bound estimates for each of the quotients here as follows
$$
\align
0 < \frac{\vartheta_3(v_n, 2i\alpha)}{\vartheta_3(0, 2i\alpha)} 
& \le 1
\\
0 < \frac{\vartheta_3(\pi t', 2i\alpha)}{\vartheta_3(\pi s, 2i\alpha)}
& \le 
\frac{\vartheta_3(0, 2i\alpha)}{\vartheta_3(\tfrac\pi2, 2i\alpha)}
\\
0 < \frac{\vartheta_3(\pi t', i\tau_\alpha)}{\vartheta_3(\pi s, i\tau_\alpha)} 
\le \frac{\vartheta_3(0, i\tau_\alpha)}{\vartheta_3(\tfrac\pi2, i\tau_\alpha)}
&\le \frac{\vartheta_3(0, 4i)}{\vartheta_3(\tfrac\pi2, 4i)}
 \le 1.00001.
\endalign
$$
The first two inequalities are obvious.  The third inequality follows from the fact 
that $\vartheta_3(0,ix)$ is a decreasing function of $x$ and 
$\vartheta_3(\tfrac\pi2, ix)$ is an increasing function for, say, $x>0.2$, and since 
$\tau_\alpha \ge 4$ for all $\alpha>0$.  The latter assertion can be seen from the 
expansion
$$
\vartheta_3(\tfrac\pi2, ix) = 1 - 2\left( 
[e^{-\pi x} - e^{-\pi 2^2 x}] + [e^{-\pi 3^2 x} - e^{-\pi 4^2 x}] + \dots \right)
$$
and noting that each of the bracketed functions is a decreasing function for $x>0.2$
as can be verified directly.  One thus gets
$$
\left| \frac{\psi_n(t)}{\psi_0(s)} \right|
\ \le \
1.00001 \frac{\vartheta_3(0, 2i\alpha)}{\vartheta_3(\tfrac\pi2, 2i\alpha)} 
\frac{|Q_n(t)|}{Q_0(s)}.
$$
Using (F1) and (F2) one gets
$$
|Q_n(t)| \le 1.09 \frac{\vartheta_2(0, 2i\alpha)}{\vartheta_3(0, 2i\alpha)},
\qquad
Q_0(t) > 0.91 \frac{\vartheta_2(0, 2i\alpha)}{\vartheta_3(0, 2i\alpha)}.
$$
This proves the following result.

\proclaim{Lemma 4.3} For all integers $n$, one has
$$
\|\psi_0^{-1}\psi_n\| \le 1.19782 
\frac{\vartheta_3(0, 2i\alpha)}{\vartheta_3(\tfrac\pi2, 2i\alpha)}.
$$
\endproclaim

Therefore one gets
$$
\|\psi_0^{-1} \sqrt{2\alpha} \Dpinner h h - I \| 
=
\left\| \sum_{n\not=0} e^{-\pi\alpha\beta^2 n^2/2} \psi_0^{-1}\psi_n V_2^n \right\|
\le
1.19782
\frac{\vartheta_3(0, 2i\alpha)}{\vartheta_3(\tfrac\pi2, 2i\alpha)}
\cdot [\vartheta_3(0, \tfrac i2 \alpha\beta^2) - 1].
$$
The right side of this inequality, as a function of $\alpha>0$ (recall $\beta^2 =
4(\alpha^2+1)$), is less than 1 for $\alpha > 0.2568$, i.e. for $0 < \theta < 0.2345$
or $\beta>2.065$, and therefore we conclude that for this $\theta$-range, the element 
$\Dpinner h h$ is invertible.  This completes the proof of $(ii)$ in Theorem 3.1 
and, in particular, yields Theorem 1.3.

\bigpagebreak

\remark{Remark}
One can in fact show, with a little extra effort, that for $0 < \theta \le 0.2427$ if 
one calculates the norms $\|\psi_0^{-1}\psi_n\|$ more precisely (for the first few 
$n$'s, say) that the the element $\Dpinner h h$ is invertible (since the norm of the
preceding sum is still less than 1).
However, if one does the same for $ \theta > 0.2451$ the norm can be greater 
than 1.  (Of course, this does not mean that $\Dpinner h h$ is singular, but that
the methods here would have to be modified.)
\endremark

\bigpagebreak

We conclude by obtaining bounds on the function $\psi_0$ and describe its asymptotic 
behaviour in $\alpha$, a result that will be used in future work ([\SWd]).

\proclaim{Proposition 4.4} For $\alpha > \tfrac14$ and all real $t$ one has
$4 < \psi_0(t) < 18.$  In addition, $\psi_0(t)\to 8$ uniformly in $t$ as $\alpha\to\infty$.
\endproclaim
\demo{Proof} From above we had
$$
\psi_0(t) = 4e^{\pi\alpha/2}
\vartheta_3(0, 2i\alpha)\, \vartheta_3(\pi t, 2i\alpha)\,
\vartheta_3(\pi t, i\tau_\alpha) \cdot Q_0(t).
$$
Let
$$
E(t) = \frac{\vartheta_2(\pi t, i\tau_\alpha)}{\vartheta_3(\pi t, i\tau_\alpha)}
$$
which goes to zero uniformly in $t$ as $\alpha\to\infty$.  By (F1) one has
$|E(t)|<0.09$.  By (F2)
one has
$$
[1-E(t)] \frac{\vartheta_2(0, 2i\alpha)}{\vartheta_3(0, 2i\alpha)}
\ < \ Q_0(t) \ < \
[1+E(t)]  \frac{\vartheta_2(0, 2i\alpha)}{\vartheta_3(0, 2i\alpha)}
$$
hence
$$
4[1-E(t)]\, R(t) \ < \ \psi_0(t) \ < \ 4[1+E(t)]\, R(t)
\tag{$\dagger$}
$$
where
$$
R(t) = e^{\pi\alpha/2}
\vartheta_2(0, 2i\alpha)\, \vartheta_3(\pi t, 2i\alpha)\,
\vartheta_3(\pi t, i\tau_\alpha).
$$
From the inequalities
$$
\align
2\ <\ e^{\pi\alpha/2} \vartheta_2(0, 2i\alpha) \ &< \ 2 \vartheta_3(0, 2i\alpha)
\\
\vartheta_3(\tfrac\pi2, 2i\alpha)\  \le \ \vartheta_3(\pi t, 2i\alpha)
& \le \vartheta_3(0, 2i\alpha)
\\
\vartheta_3(\tfrac\pi2, i\tau_\alpha) \le \vartheta_3(\pi t, i\tau_\alpha)
&\le \vartheta_3(0, i\tau_\alpha)
\endalign
$$
one obtains
$$
2 \vartheta_3(\tfrac\pi2, 2i\alpha) \vartheta_3(\tfrac\pi2, i\tau_\alpha)
\ < \ R(t) \ < \ 2 \vartheta_3(0, 2i\alpha)^2 \vartheta_3(0, i\tau_\alpha).
\tag*
$$
%%%% As $x\mapsto \vartheta_3(\tfrac\pi2, ix)$ is an increasing function, the left side
Since it is not hard to see that the left side of (*) is an increasing function of
$\alpha\ge\tfrac14$, it is greater that its value at $\alpha=\tfrac14$ which, when 
computed (and using $\tau_\alpha\ge 4$), gives a value greater than $1.17$.
The right side of (*) is a decreasing function of $\alpha$ (as is each factor), so 
is at most equal to its value when $\alpha=\tfrac14$, which is less than $4.03$.
Therefore, $1.17<R(t)<4.03$ for all $t$ and $\alpha\ge\tfrac14$, which yields the 
claimed inequality of the lemma for $\psi_0$.  
Finally, note that the left and right sides of (*) converge
to $2$ as $\alpha\to\infty$, so that $R(t)\to2$ uniformly in $t$.
Therefore from ($\dagger$) one has $|\psi_0(t)-4R(t)| \le |E(t)|R(t)$ and
thereby obtains the claimed asymptotic behaviour of $\psi_0$.  \qed
\enddemo

\bigpagebreak

%%%%%%%%%%%%%%%%%%%%%%%%%%%%%%%%%%%%    REFERENCES      %%%%%%%%%%%
%%%%%%%%%%%%%%%%%%%%%%%%%%%%%%%%%%%%    REFERENCES      %%%%%%%%%%%
%%%%%%%%%%%%%%%%%%%%%%%%%%%%%%%%%%%%    REFERENCES      %%%%%%%%%%%
%%%%%%%%%%%%%%%%%%%%%%%%%%%%%%%%%%%%    REFERENCES      %%%%%%%%%%%
%%%%%%%%%%%%%%%%%%%%%%%%%%%%%%%%%%%%    REFERENCES      %%%%%%%%%%%
%%%%%%%%%%%%%%%%%%%%%%%%%%%%%%%%%%%%    REFERENCES      %%%%%%%%%%%

\enddocument